\documentclass[12pt]{article}
\usepackage{amssymb}
\usepackage{graphicx}
\usepackage{amsfonts}
\usepackage{latexsym}
\usepackage{amsthm}
\usepackage{amsmath}

\usepackage{amssymb, amscd}

\usepackage{url}
\usepackage[T1]{fontenc}

\newtheorem{problem}{Problem}
\newtheorem{theo}[problem]{Theorem}
\newtheorem{rem}[problem]{Remark}

\newtheorem{defin}[problem]{Definition}
\newtheorem{prop}[problem]{Proposition}
\newtheorem{cor}[problem]{Corollary}
\newtheorem{lema}[problem]{Lemma}
\newtheorem{exam}[problem]{Example}

\newtheorem{conj}[problem]{Conjecture}

\begin{document}
\date{May 30, 2013}
\title{ Measurable patterns, necklaces,\\ and sets indiscernible by measure}

\author{{\Large Sini\v sa Vre\' cica}\thanks{Supported by the Ministry
of Education and Science of Serbia, Grant 174034.}\\ {\small Faculty of Mathematics}\\[-2mm] {\small University of Belgrade}
  \and {\Large Rade  \v Zivaljevi\' c}\thanks{Supported by the Ministry
  of Education and Science of Serbia, Grants 174020  and  174034.}\\ {\small Mathematical Institute}\\[-2mm] {\small SASA,
  Belgrade}\\[-2mm]}

\maketitle
\begin{abstract}
In some recent papers the classical `splitting necklace theorem'
is linked in an interesting way with a geometric `pattern
avoidance problem', see Alon et al.\ (Proc.\ Amer.\ Math.\ Soc.,
2009), Grytczuk and Lubawski (arXiv:1209.1809 [math.CO]), and
Laso\'{n} (arXiv:1304.5390v1 [math.CO]). Following these authors
we explore the topological constraints on the existence of a
(relaxed) measurable coloring of $\mathbb{R}^d$ such that any two
distinct, non-degenerate cubes (parallelepipeds) are measure
discernible. For example, motivated by a conjecture of Laso\'{n},
we show that for every collection $\mu_1,...,\mu_{2d-1}$ of $2d-1$
continuous finite measures on $\mathbb{R}^d$, there exist two
nontrivial axis-aligned $d$-dimensional cuboids (rectangular
parallelepipeds) $C_1$ and $C_2$ such that $\mu_i(C_1)=\mu_i(C_2)$
for each $i\in\{1,...,2d-1\}$. We also show by examples that the
bound $2d-1$ cannot be improved in general. These results are
steps in the direction of studying general topological
obstructions for the existence of non-repetitive colorings of
measurable spaces.
\end{abstract}

\section{Introduction}
\label{sec:intro}

The following definition explains in what sense two objects
(measurable sets) can be {\em measure discernible} (or {\em
indiscernible}).

\begin{defin}\label{def:indiscernible}
Let $(X, \mathcal{B}, \mu_1,\ldots ,\mu_d)$ be a measure space
with a collection $\mu =(\mu_1,\ldots,\mu_d)$ of measures. We say
that two measurable sets $A, B\in \mathcal{B}$ are
$\mu$-indiscernible or {\em measure indiscernible} if $\mu(A) =
\mu(B)$ as vectors in $\mathbb{R}^d$ or equivalently if,
\begin{equation}\label{eqn:indiscernible}
\mu_j(A) = \mu_j(B) \qquad \mbox{ {\rm for each} } j=1,\ldots, d.
\end{equation}
In the opposite case, i.e.\ if at least one of equalities in
(\ref{eqn:indiscernible}) does not hold, the sets $A$ and $B$ are
measure discernible.
\end{defin}

There are many interesting combinatorial geometric results which
claim the existence of measure indiscernible partitions of the
ambient space $X$. The classical `Ham Sandwich Theorem' is a
result of this type. Indeed, if $X = \mathbb{R}^d$ then it claims
the existence of two $\mu$-indiscernible half-spaces which have a
common boundary hyperplane. Much more recent is the result of
Hubard and Aronov \cite{Hub-Aro}, Karasev \cite{Karasev}, and
Sober\'{o}n \cite{Soberon}, who showed that for a given collection
of $d$ continuous measures $\mu_1,\ldots,\mu_d$, defined on
$\mathbb{R}^d$, and an integer $k\geq 2$, there exists a partition
of $\mathbb{R}^d$ into $k$ convex sets which are
$\mu$-indiscernible.

\medskip
Some questions (and results) about indiscernible partitions are
better known as problems about {\em fair division}, {\em consensus
partitions}, {\em envy free divisions},  or simply as {\em
equipartitions} of measures, \cite{M, Long, MVZ, Z04, Z08}. One of
the best known results of this type is the {\em `splitting
necklace theorem'} of Alon \cite{Alo87, alon-constructive} which
says that each necklace with $k\cdot a_i$ beads of color
$i=1,\ldots, n$ can be fairly divided between $k$ \emph{thieves}
by at most $n(k-1)$ cuts. Alon deduced this result from the fact
that such a division is possible also in the case of a continuous
necklace $[0,1]$ where beads of given color are interpreted as
measurable sets $A_i\subset [0,1]$ (or more generally as
continuous measures $\mu_i$).

\medskip
Some `pattern avoidance problems' \cite{Currie} also appear to be
directly related to questions about measure indiscernible sets,
however until recently \cite{AGLM, GL, Las}  these areas seem to
have had completely independent development. For illustration,
Erd\"{o}s \cite{Erdos} asked whether there is a $4$-coloring of
the integers such that each two adjacent intervals are (in our
terminology) `color discernible', meaning that they remain
different even after some permutation of their elements.
Continuous (measure theoretic) analogues of these questions were
formulated and studied in \cite{AGLM}, \cite{GL} and \cite{Las}.

\medskip
The paper \cite{AGLM} establishes an interesting link between the
pattern avoidance problem of Erd\"{o}s and the splitting necklaces
problem and focuses on the question whether the number of cuts can
be reduced for some subinterval of a line measurably colored by a
prescribed number of colors. For example they showed that there
exists a measurable $4$-coloring of the real line such that two
adjacent intervals are always color discernible.

\medskip
Papers \cite{GL} and \cite{Las} continued this research,
connecting the higher dimensional pattern avoidance problem with
the higher dimensional extensions of the splitting necklace
theorem \cite{Long-Ziv}. In particular the results and conjectures
of Laso\'{n}  \cite{Las} are our immediate motivation for
exploring these and other aspects of measurable colorings of
Euclidean spaces.



\subsection{Our paper}

Our objective is to identify and explore the topological
constraints for the existence of `non-repetitive' or `pattern
avoiding' colorings (measures) of $\mathbb{R}^d$ which are not
necessarily measurable partitions
(Definition~\ref{def:measurable-colorings}). For a given family of
measurable sets in $\mathbb{R}^d$ we introduce the
`pattern-avoiding number' $\nu(\mathcal{F})$ and the `relaxed
pattern-avoiding number' $\nu_{rel}(\mathcal{F})$
(Definition~\ref{def:nu-invariant}) which detect the critical
number of colors when color repetitions in $\mathcal{F}$ are
always present.

\medskip
Some results and conjectures of Laso\'{n} \cite{Las} (see
Section~\ref{sec:Lason} for an outline) are naturally interpreted
as results about the invariant $\nu(\mathcal{F})$. Our focus is on
the closely related invariant $\nu_{rel}(\mathcal{F})$ which is
easier to handle so we are able to provide much more precise
information, including some exact calculations.

\medskip
Following  \cite{AGLM, GL, Las} we put some emphasis on the class
$\mathcal{C}_d$ of $d$-cubes and the class $\mathcal{P}_d$ od
$d$-cuboids (rectangular parallelepipeds) in $\mathbb{R}^d$. Our
first exact evaluation (Theorem~\ref{thm:prva},
Examples~\ref{exam:prvi-primer}, \ref{exam:drugi-primer}, and
\ref{exam:treci-primer}) shows that $$\nu_{rel}(\mathcal{C}_d)=
d+1 \qquad {\rm and} \qquad  \nu_{rel}(\mathcal{P}_d)= 2d.$$

The class $\mathcal{P}_d$ is somewhat special in the sense that it
is invariant with respect to a very large group of
auto-homeomorphisms of $\mathbb{R}^d$
(Remark~\ref{rem:justify-cubes}).  As a consequence one can
calculate the generalized $\nu$-invariant (in the sense of
Remark~\ref{rem:general-nu}) for the class $\mathcal{P}_d$ not
only for signed measures but for some other classes  including the
positive and probability measures on $\mathbb{R}^d$.

\medskip
In other directions we show that Theorem~\ref{thm:prva} admits
several extensions of different nature. In Theorem~\ref{thm:druga}
we prove, by using more powerful topological tools, that one can
often guarantee the existence of an arbitrarily large finite
family of measure indiscernible cubes (cuboids). The same method
yields an even stronger result involving families obtained by more
general Lie group actions (Theorem~\ref{thm:Lie-group}).

\medskip
In the special case of the Lie group $G_{DL}$, generated by
positive, axis-aligned dilatations and translations in
$\mathbb{R}^d$, we calculate (Theorem~\ref{thm:opposite}) the
$\nu_{rel}$-invariant of the associated families of centrally
symmetric convex bodies.

\medskip
We study also other aspects of Theorem~\ref{thm:prva} and show for
example (Section~\ref{sec:disjoint}) that in same instances of the
problem one can guarantee the existence of {\em disjoint} measure
indiscernible cubes.

\section{Non-repetitive colorings of $\mathbb{R}^d$}

If not specified otherwise, all measures are signed, finite Borel
measures on $\mathbb{R}^d$ which are absolutely continuous with
respect to the Lebesgue measure $dm$.

\begin{defin}\label{def:measurable-colorings}
A measurable $k$-coloring of $\mathbb{R}^d$ is a partition
$\mathbb{R}^d = A_1\cup\ldots\cup A_k$ of the ambient $d$-space
into $k$-measurable sets. A {\em relaxed} $k$-coloring of
$\mathbb{R}^d$ is a collection $\mu = (\mu_1,\ldots,\mu_k)$ of
continuous measures, $d\mu_i = f_i\, dm$, where $f_i(x)$ is the
`intensity' of color $i$ at $x$, and (more importantly) $\mu_i(A)$
is the total amount of color $i$ used for coloring of the
measurable set $A$.
\end{defin}

The following definition explains in what sense a (relaxed)
measurable $k$-coloring of $\mathbb{R}^d$ may be {\em
non-repetitive} (pattern avoiding). The definition is formulated
in the language of measures (relaxed colorings) but we tacitly use
it also for strict colorings (partitions).

\begin{defin}\label{def:pattern-avoiding}
Let $\mathcal{F}$ be a family of Lebesgue measurable sets in
$\mathbb{R}^d$, such as the family of all axis-aligned cubes\,
$\mathcal{C}_d$ or the family $\mathcal{P}_d$ of all axis-aligned
cuboids ($d$-parallelepipeds). We say that a (relaxed)
$k$-coloring $\mu = (\mu_1, \ldots, \mu_k)$ of\, $\mathbb{R}^d$ is
$\mathcal{F}$-non-repetitive (cube non-repetitive, cuboid
non-repetitive) if each two distinct elements $A, B\in
\mathcal{F}$ are $\mu$-discernible
(Definition~\ref{def:indiscernible}) i.e.\ if $\mu_i(A)\neq
\mu_i(B)$ for at least one of the indices $i\in [k]$.
\end{defin}

\begin{defin}\label{def:nu-invariant}
Given a family $\mathcal{F}$ of Lebesgue measurable sets in
$\mathbb{R}^d$ we define the corresponding (measure)
`pattern-avoiding number' of $\mathcal{F}$ as the number,
\begin{equation}\label{eqn:nu-inv}
\nu(\mathcal{F})  = \mbox{\rm Inf}\{k \in \mathbb{N} \mid
\exists\, \mathcal{F}\mbox{-non-repetitive }\, k \mbox{-coloring
of }\, \mathbb{R}^d\}.
\end{equation}
Similarly, if we alow relaxed colorings we have the corresponding
`relaxed pattern avoiding number' $\nu_{rel}(\mathcal{F})$ defined
as the minimum (infimum) of all $k$ such that there exists a
relaxed $\mathcal{F}$-non-repetitive coloring of $\mathbb{R}^d$.
\end{defin}

\begin{rem}\label{rem:general-nu}{\rm
Perhaps a more systematic approach would involve more general
invariants $\nu(\mathcal{F}, \mathcal{C})$ where, aside from the
family $\mathcal{F}$ of measurable sets, one also specifies in
advance the family $\mathcal{C}$ of admissible colorings
(measures). Here we deal mainly with `relaxed colorings' and the
corresponding invariant $\nu_{rel}(\mathcal{F})$, where
$\mathcal{C}$ is the class of all signed, continuous, finite Borel
measures. The invariant $\nu(\mathcal{F})$ is recovered if
$\mathcal{C}$ is the family of all colorings with disjoint
measurable set. Other cases of interest would include positive
(probabilistic) measures, measures satisfying a condition on their
support, etc.}
\end{rem}

\subsection{Some known results about $\nu(\mathcal{F})$}
\label{sec:Lason}

Continuing the research from \cite{AGLM} and \cite{GL}, and in
particular improving over some bounds established in \cite{GL},
Laso\'{n} in \cite{Las} described a method of constructing
measurable $k$-colorings of $\mathbb{R}^d$ which are cube (or
cuboid) non-repetitive.

\begin{theo}\label{thm:Lason-3.6}{\rm \cite[Theorem~3.6]{Las}}
For every $d \geq 1$ there exists a measurable $(2d + 3)$-coloring
of $\mathbb{R}^d$ such that no two nontrivial axis-aligned
d-dimensional cubes have the same measure of each color. In other
words there exists a measurable $(2d+3)$-coloring (partition) of
$\mathbb{R}^d$ which is $\mathcal{C}_d$-non-repetitive.
\end{theo}

\begin{theo}\label{thm:Lason-3.8}{\rm \cite[Theorem~3.8]{Las}}
There exists a measurable $(4d+1)$-coloring (partition) of
$\mathbb{R}^d$ which is $\mathcal{P}_d$-non-repetitive. In other
words there exists a measurable $(4d + 1)$-coloring of
$\mathbb{R}^d$ such that no two nontrivial axis-aligned
$d$-dimensional cuboids have the same measure of each color.
\end{theo}

It is natural to ask whether the bounds $2d+3$ and $4d+1$ in
Theorems~\ref{thm:Lason-3.6} and \ref{thm:Lason-3.8} are the best
possible so Laso\'{n} formulated also the following conjecture.

\begin{conj}\label{conj:Lason}{\rm \cite[Conjecture~3.7]{Las}}
For every measurable $(2d + 2)$-coloring of $\mathbb{R}^d$ there
exist two non-degenerate axis-aligned $d$-dimensional cubes which
have the same measure of each color. In other words each $(2d +
2)$-coloring of $\mathbb{R}^d$ is `pattern-repetitive' in the
sense that there always exist two distinct cubes that are measure
indiscernible.
\end{conj}

\subsection{Repetitive relaxed colorings of $\mathbb{R}^d$}

Here we address the question of the existence of cube (cuboid)
non-repetitive patterns in the class of {\em relaxed} measurable
colorings of $\mathbb{R}^d$
(Definition~\ref{def:measurable-colorings}). In other words we
consider exactly the same questions addresses by
Theorems~\ref{thm:Lason-3.6} and \ref{thm:Lason-3.8} and
Conjecture~\ref{conj:Lason} but we allow more general colorings
provided by continuous measures which do not necessarily
correspond to measurable partitions.

\medskip
Aside from proving the counterparts of
Theorems~\ref{thm:Lason-3.6} and \ref{thm:Lason-3.8} we also
provide examples showing that the bounds are the best possible in
this case.

\begin{theo}\label{thm:prva}
For every collection $\mu_1,...,\mu_d$ of $d$ continuous finite
measures on $\mathbb{R}^d$, there are two nontrivial axis-aligned
$d$-dimensional cubes $C_1$ and $C_2$ such that
$\mu_i(C_1)=\mu_i(C_2)$ for all $i=1,...,d$.

For every collection $\mu_1,...,\mu_{2d-1}$ of $2d-1$ continuous
finite measures on $\mathbb{R}^d$, there are two nontrivial
axis-aligned $d$-dimensional cuboids (rectangular parallelepiped)
$C_1$ and $C_2$ such that $\mu_i(C_1)=\mu_i(C_2)$ for all
$i=1,...,2d-1$.
\end{theo}

\medskip\noindent
{\bf Proof:} Each nontrivial axis-aligned cube in $\mathbb{R}^d$
is uniquely determined by its vertex $a=(a_1,...,a_d)$ with
smallest coordinates and with the length $l$ of its edge. So, the
space of all such cubes is homeomorphic to $\mathbb{R}^d\times
(0,\infty )$.

The configuration space of all pairs of distinct, axis-aligned
cubes in $\mathbb{R}^d$ can be described as $(\mathbb{R}^d\times
(0,\infty ))^2\setminus \Delta$, where $\Delta$ is the diagonal in
the product space. This space is obviously
$\mathbb{Z}/2$-equivariantly homotopy equivalent to the sphere
$S^d$. (The antipodal action on the configuration space is the
obvious one.)

Let us consider the $\mathbb{Z}/2$-equivariant mapping $F :
(\mathbb{R}^d\times (0,\infty ))^2\setminus \Delta \rightarrow
\mathbb{R}^d$ given by
$$F((a,l_1),(c,l_2))=(\mu_1(a,l_1)-\mu_1(c,l_2),...,
\mu_d(a,l_1)-\mu_d(c,l_2)).$$

If there are no measure indiscernible cubes, this mapping would
miss the origin. This would lead to an antipodal map from $S^d$ to
$S^{d-1}$, which is a contradiction establishing the `cube case'
of the theorem.

\medskip
Each nontrivial axis-aligned cuboid in $\mathbb{R}^d$ is uniquely
determined by its vertex $a=(a_1,...,a_d)$ with smallest
coordinates and with its vertex $b=(b_1,...,b_d)$ with biggest
coordinates. Here $a_i<b_i$ for each $i=1,...,d$. So, for every
$i$, $(a_i,b_i)$ always belongs to the open half-plane $P$ (above
the line $y=x$). Therefore, the space of all such cuboids is
homeomorphic to $\mathcal{P}_d\approx (\mathbb{R}^2)^d$.

As a consequence, the configuration space of all pairs of two
different nontrivial axis-aligned cuboids in $\mathbb{R}^d$ can be
described as $(\mathcal{P}_d)^2\setminus \Delta$, where $\Delta$
is the diagonal in the product space. This space is obviously
$\mathbb{Z}/2$-equivariantly homotopy equivalent to the sphere
$S^{2d-1}$.

Let us consider the $\mathbb{Z}/2$-equivariant mapping $G :
(\mathcal{P}_d)^2\setminus \Delta \rightarrow \mathbb{R}^{2d-1}$
given by
\begin{equation}\label{eqn:antipod}
F((a,b),(c,d))=(\mu_1(a,b)-\mu_1(c,d),...,
\mu_{2d-1}(a,b)-\mu_{2d-1}(c,d)).
\end{equation}

If there are no pairs of measure indiscernible cuboids, the map
described by (\ref{eqn:antipod}) would miss the origin. This would
imply the existence of an antipodal map from $S^{2d-1}$ to
$S^{2d-2}$, which leads to the desired contradiction. \hfill
$\square$
\bigskip

We complete this section by providing examples showing that the
estimates obtained in Theorem~\ref{thm:prva} are the best
possible. We describe the densities of continuous {\em signed}
measures on $\mathbb{R}^d$ which restricted on $I^d = (0,1)^d$
yield (after normalizing) probability measures on the open
$d$-cube $I^d$. In light of Remark~\ref{rem:justify-cubes} these
measures can be pulled back to $\mathbb{R}^d$ to yield desired
probability measures on $\mathbb{R}^d$.

\begin{exam}\label{exam:prvi-primer}
{\rm For illustration we initially treat the case $d=1$, and give
two measures on the segment $(0,1)$ such that no two different
intervals in $(0,1)$ contain the same amount of both measures.
Notice that in the case $d=1$ both parts of Theorem~\ref{thm:prva}
reduce to the same statement.

The measures $\mu_1$ and $\mu_2$ are described by their density
functions $\varphi_1(x)=1-x$ and $\varphi_2(x)=x$. Any two
intervals containing the same amount of both measures would be of
the same length, since the density of the measure $\mu_1+\mu_2$ is
constant. But, if they are different, one of them would be more
"on the left" and that one would contain greater amount of the
measure $\mu_1$ and smaller amount of the measure $\mu_2$ than the
other interval.

Notice that this example is also related to the conjecture 3.4 in
\cite{Las}. Namely, it is conjectured there that for any partition
of $\mathbb{R}^d$ in $k$ measurable sets, there is an axis-aligned
cube which has a fair $q$-splitting using at most $k(q-1)-d-1$
axis-aligned hyperplane cuts. This example shows that in the case
of measures (and not partitions) for $k=q=2$ and $d=1$, one cut is
not enough, while $k(q-1)-d-1=0$.}
\end{exam}

\begin{exam}\label{exam:drugi-primer}
{\rm Let us now take care of the general case. We construct first
a collection of $d+1$ measures on $(0,1)^d$ such that no two
different cubes contain the same amount of every measure. Let
these measures be given by their density functions
$\varphi_1,...,\varphi_{d+1} : (0,1)^d\rightarrow \mathbb{R}$,
$\varphi_1(x_1,...,x_d)=x_1x_2\cdots x_d$,
$\varphi_2(x_1,...,x_d)=(1-x_1)x_2\cdots x_d$,
$\varphi_3(x_1,...,x_d)=(1-x_2)x_3\cdots x_d$,
$\varphi_4(x_1,...,x_d)=(1-x_3)x_4\cdots x_d$,...,
$\varphi_{d+1}(x_1,...,x_d)=1-x_d$.

It is easy to verify that the length of the edges of two cubes
containing the same amount of every measure have to be equal, and
than that the coordinates of the vertex with the smallest
coordinates should be also equal for these two cubes. (We first
notice that for the coordinate $x_d$, then $x_{d-1}$ etc.}
\end{exam}

\begin{exam}\label{exam:treci-primer}
{\rm Let us describe $2d$ measures on $\mathbb{R}^d$ by their
density functions $\varphi_0(x)=1$, $\varphi_i(x)=x_i$ for $i\in
\{1,...,d\}$, and $\varphi_{d+i}(x)=x_i^3$ for $i\in
\{1,...,d-1\}$.

Let us denote by $[a,b]$ and $[c,d]$ two cuboids containing the
same amount of every measure. Here $a=(a_1,...,a_d)$ and
$c=(c_1,...,c_d)$ are vertices of these two cuboids with smallest
coordinates and $b$ and $d$ the vertices with biggest coordinates.

The requirement that these cuboids contain the same amount of
every measure provide us with the following equalities.

Measure $\mu_0$ gives us the equality
$\prod_{i=1}^d(b_i-a_i)=\prod_{i=1}^d(d_i-c_i)$. Measures $\mu_j$
for $j\in \{1,...,d\}$ give us the equalities $(a_j+b_j)
\prod_{i=1}^d(b_i-a_i)=(c_j+d_j)\prod_{i=1}^d(d_i-c_i)$. Together
with the first equality these give us the equalities
$a_j+b_j=c_j+d_j$, for all $j\in \{1,...,d\}$. Measures
$\mu_{d+j}$ for $j\in \{1,...,d-1\}$ give us the equalities
$(a_j^2+b_j^2)(a_j+b_j)
\prod_{i=1}^d(b_i-a_i)=(c_j^2+d_j^2)(c_j+d_j)\prod_{i=1}^d(d_i-c_i)$.
Together with previous equalities these give us
$a_j^2+b_j^2=c_j^2+d_j^2$ for $j\in \{1,...,d-1\}$. Furthermore,
these equalities together with the equalities $a_j+b_j=c_j+d_j$
give us directly $a_j=c_j$ and $b_j=d_j$ for $j\in \{1,...,d-1\}$.
Then, from the first equality we have $b_d-a_d=d_d-c_d$, and from
another one $a_d+b_d=c_d+d_d$. This gives us $a_d=c_d$ and
$b_d=d_d$, and so these two cuboids coincide. This means that two
different cuboids could not contain the same amount of every of
the described $2d$ measures.

\begin{cor}\label{cor:best-possible}
If $\mathcal{C}_d$ and $\mathcal{P}_d$ are the families of all
cubes (respectively cuboids) in the $d$-dimensional space
$\mathbb{R}^d$ then,
\begin{equation}\label{eqn:tacna-nu-vrednost}
\nu_{rel}(\mathcal{C}_d) = d+1 \qquad and \qquad
\nu_{rel}(\mathcal{P}_d) = 2d.
\end{equation}
\end{cor}

Let us compare the results obtained in this section with the
results of \cite{Las} (see also Section~\ref{sec:Lason}), dealing
(instead of measures) with the partitions of $\mathbb{R}^d$ in
disjoint measurable subsets. Laso\'{n} proved that there exists a
partition of $\mathbb{R}^d$ in $2d+3$ disjoint measurable sets
such that no two different nontrivial axis-aligned cubes contain
the same amount of every of these sets. Also, it is proved that
there exists a partition of $\mathbb{R}^d$ in $4d+1$ disjoint
measurable sets such that no two different nontrivial axis-aligned
cuboids contain the same amount of every of these sets. It is
conjectured (Conjecture~\ref{conj:Lason} in
Section~\ref{sec:Lason}) that these estimates are the best
possible.

We work with measures and prove the corresponding results with
$2d+3$ being replaced by $d+1$, and with $4d+1$ being replaced by
$2d$, and show that these results are the best possible in this
case.}
\end{exam}

\section{A generalization}

\begin{defin}\label{def:configuration-space} For each topological
space $X$, the associated {\em configuration space} $F(X,n)$ of
all $n$-tuples of labelled points in $X$ is the space,
$$
F(X,n) := \{x\in X^n \mid x_i\neq x_j \mbox{ {\rm for each} }
i\neq j \}.
$$
\end{defin}
The obvious action of the symmetric group $S_n$ on $X^n$,
restricts to a free action on the associated configuration space
$F(X,n)$.

\medskip

As shown by examples in the previous section,
Theorem~\ref{thm:prva} is optimal as far as the number of measures
is concerned. However, we show that it can be considerably
improved in a different direction. Indeed, it turns out that
instead of two cubes (cuboids) we can prove the existence of a
finite family of cubes (cuboids) of any size which are
$\mu$-indiscernible in the sense of
Definition~\ref{def:indiscernible}.

More explicitly, we extend the results from the previous section
to the case of $n$ cubes (cuboids) in $\mathbb{R}^d$ where $n=p^k$
is a power of a prime $p$. Instead of the Borsuk-Ulam theorem,
used in the proof of Theorem~\ref{thm:prva}, we apply the
following result, see \cite{Vas, Hung, Grom, Karasev, Hub-Aro,
Bla-Zie}.

\begin{theo}\label{thm:BU-general}
Suppose that $n=p^k$ is a power of a prime $p\geq 2$ and let
$m\geq 2$. Let $W_n$ be the $(n-1)$-dimensional, real
representation of the symmetric group $S_n$ which arises as the
orthogonal complement of the diagonal in the permutation
representation $\mathbb{R}^n$.  Then each equivariant map
\begin{equation}\label{eqn:eq-map}
\Phi : F(\mathbb{R}^m, n) \rightarrow W_n^{\oplus (m-1)}
\end{equation}
must have a zero.
\end{theo}

\begin{theo}\label{thm:druga}
For each collection $\mu_1,...,\mu_d$ of $d$ continuous, finite
measures on $\mathbb{R}^d$ and any natural number $n$, there
exists a collection of $n$ nontrivial, axis-aligned
$d$-dimensional cubes $C_1, C_2,\ldots, C_n$ which are
$\mu$-indiscernible in the sense that $\mu_i(C_j)=\mu_i(C_k)$ for
all $i=1,...,d$ and $j,k\in\{1,\ldots,n\}$.

For each collection $\mu_1,...,\mu_{2d-1}$ of $2d-1$ continuous
finite measures on $\mathbb{R}^d$, there exists a collection of
$n$ nontrivial, axis-aligned $d$-dimensional cuboids (rectangular
parallelepipeds) $C_1, C_2,\ldots, C_n$ which are
$\mu$-indiscernible in the sense that $\mu_i(C_j)=\mu_i(C_k)$ for
all $i=1,...,2d-1$ and $j,k\in\{1,\ldots,n\}$.
\end{theo}

\medskip\noindent
{\bf Proof:} We outline the proof of the second statement.

Without loss of generality we can assume that $n = p^k$ is a power
of a prime number (say $n=p$ for a prime $p$). As already observed
in the proof of Theorem~\ref{thm:prva}, the variety
$\mathcal{P}_d$ of all cuboids in $\mathbb{R}^d$ is homeomorphic
to $\mathbb{R}^{2d}$. Given a collection $C = (C_1,\ldots, C_n)\in
F(\mathcal{P}_d,n)$ of pairwise distinct cuboids and a measure
$\mu_i$ let $\mu_i(C) := (\mu_i(C_1),\ldots, \mu_i(C_n))\in
\mathbb{R}^n$. Obviously the cuboids $\{C_j\}_{j=1}^n$ are
$\mu_i$-indiscernible if and only if $\pi(\mu_i(C))=0$ where $\pi
: \mathbb{R}^n \rightarrow W_n$ is the natural projection.

Let $\phi_i : F(\mathcal{P}_d, n)\rightarrow W_n$ be the map
defined by $\phi_i(C) = \pi(\mu_i(C))$. Let $\Phi :
F(\mathcal{P}_d, n)\rightarrow (W_n)^{\oplus (2d-1)}$ be the
associated map where $\Phi(C) = (\phi_1(C), \ldots,
\phi_{2d-1}(C))$. Since $\mathcal{P}_d \cong \mathbb{R}^{2d}$ it
follows from Theorem~\ref{thm:BU-general} that for some $C\in
F(\mathcal{P}_d, n),\, \Phi(C)=0$ which completes the proof of the
theorem. \hfill $\square$

\subsection{The case of pairwise disjoint cubes and cuboids}
\label{sec:disjoint}

A natural question is whether one can strengthen
Theorems~\ref{thm:prva} and \ref{thm:druga} by claiming the
existence of {\em pairwise disjoint} cuboids (cubes) which are
$\mu$-indiscernible.

\begin{prop}\label{prop:May}
The configuration space of all ordered collections of $n\geq 2$
cuboids in $\mathbb{R}^d$ is $S_n$-equivariantly homotopy
equivalent to the configuration space $F(\mathbb{R}^d,n)$.
Moreover, the configuration space of all ordered collections of
$n\geq 2$ cubes in $\mathbb{R}^d$ of the same size is also
$S_n$-equivariantly homotopy equivalent to the configuration space
$F(\mathbb{R}^d,n)$.
\end{prop}

\medskip\noindent
{\bf Proof:} As shown by May \cite[Theorem~4.8]{May}, the
configuration space of all ordered collections of $n$ pairwise
disjoint, axis-aligned cuboids in $\mathbb{R}^d$ is
$S_n$-equivariantly homotopy equivalent to the configuration space
$F(\mathbb{R}^d, n)$. Actually May puts more emphasis in his proof
on axis-aligned cuboids with disjoint interiors but the argument
in \cite{May} can be easily modified to cover the case of disjoint
cuboids as well.

A similar result holds for disjoint cubes. Indeed, for a given
axis-aligned cuboid $C$ let $\widehat{C}$ be the
$\subseteq$-maximal cube in the set of all axis-aligned cubes
$D\subset C$ which have the same barycenter as $C$. Then the map
which sends an ordered collection $(C_1,\ldots, C_n)$ of $n$
pairwise disjoint cuboids to the corresponding collections
$(\widehat{C}_1,\ldots, \widehat{C}_n)$ of cubes is easily shown
to be a deformation retraction. Finally, if all cubes
$\widehat{C_i}$ are shrank to the cubes of the same size we obtain
a deformation retraction that establishes the second part of the
proposition. \hfill $\square$

\medskip
Proposition~\ref{prop:May} shows that we cannot improve the
`cuboid case' of Theorem~\ref{thm:druga} (by an argument based on
Theorem~\ref{thm:BU-general}) to pairwise disjoint cuboids unless
we drastically reduce the number of measures. However, the
situation with cubes is different. The following result shows that
one can always find two or more pairwise disjoint, measure
indiscernible cubes of the {\em same size} if one allows not more
than $(d-1)$ colors (i.e.\ one less than in the `cube case' of
Theorem~\ref{thm:druga}). The proof relies on
Proposition~\ref{prop:May} and follows closely the proof of
Theorem~\ref{thm:druga} so the details are omitted.

\begin{prop}\label{prop:treca}
For each collection $\mu_1,...,\mu_{d-1}$ of $(d-1)$ continuous,
finite measures on $\mathbb{R}^d$ and any natural number $n$,
there exists a collection of $n$ pairwise disjoint, axis-aligned
$d$-dimensional cubes $C_1, C_2,\ldots, C_n$ of the {\em same
size} which are $\mu$-indiscernible in the sense that
$\mu_i(C_j)=\mu_i(C_k)$ for all $i=1,...,d$ and
$j,k\in\{1,\ldots,n\}$.
\end{prop}

\begin{conj}\label{conj:best-possible}{\rm
The result from Proposition~\ref{prop:treca} is the best possible
in the following stronger sense. There exists a collection of $d$
continuous, finite, measures on $\mathbb{R}^d$ such that not only
pairs of disjoint cubes but the pairs of disjoint {\em cuboids}
are also measure discernible.}
\end{conj}

\section{Towards general non-repetitive colorings}

Here we show that there is nothing special about cubes and cuboids
and that theorems from the previous sections hold also for balls,
ellipsoids, and even more generally for measurable sets of
suitable form. Moreover, there is nothing special about choosing a
preferred position for the selected geometric shape (say axes
aligned or similar). Perhaps the most natural framework for this
problem is given by the following result which involves arbitrary
Lie group actions on $\mathbb{R}^d$.

\begin{theo}\label{thm:Lie-group}
Let $Q$ be a polytope in $\mathbb{R}^d$ (more generally a convex
body or a just a measurable set). Let $G$ be a Lie group acting on
$\mathbb{R}^d$. Let $\mathcal{F} = O_Q = \{g(Q) \mid g\in G\}$ be
the set (orbit) of all images of $C$ with respect to actions of
elements from $G$. Assume that $O_Q$ is a smooth manifold
(possibly with singularities) of geometric dimension $\nu$. Then
for each {\em relaxed measurable coloring} of $\mathbb{R}^d$ with
$\nu-1$ colors (Definition~\ref{def:measurable-colorings}) and
each integer $n\geq 2$ there exist a collection of $n$ distinct
elements in $\mathcal{F}$ which are pairwise measure indiscernible
(in the sense of Definition~\ref{def:indiscernible}). In
particular, $\nu_{rel}(\mathcal{F})\geq \nu$.
\end{theo}

\medskip\noindent
{\bf Proof:} By passing to a larger number if necessary we can
assume that $n=p^k$ is power of a prime. By assumption $O_Q$ is a
manifold of dimension $\nu$ so there is a subset $U\subset O_Q$
homeomorphic to $\mathbb{R}^\nu$. The test map for the existence
of a collection of $n$ measure indiscernible sets in $U$ is,
\[
\phi_i : F(U, n)\rightarrow W_n^{\oplus (\nu-1)}
\]
which by Theorem~\ref{thm:druga} must have a zero. \hfill
$\square$

\subsection{Remarks and examples}
Typically the orbit space $O_Q$ that appears in
Theorem~\ref{thm:Lie-group} is homeomorphic to a homogeneous
manifold $G/H$ where $H = \{g\in G \mid g(Q)=Q\}$. For example it
is well-known that the group of all isometries of a convex body
$K$ in $\mathbb{R}^d$ is a Lie subgroup $H$ of the group
$G=Isom(\mathbb{R}^d)$ of all isometries of the ambient space. In
this case Theorem~\ref{thm:Lie-group} establishes a connection
between the dimension of the isometry (symmetry) group $H$ of $K$
and the $\nu_{rel}$-invariant of the associated family
$\mathcal{F}_K$ of isometric copies of $K$ in $\mathbb{R}^d$,
\begin{equation}\label{eqn:nu-isom}
d(d+1)/2 < \mbox{\rm dim}(H) + \nu_{rel}(\mathcal{F}_K).
\end{equation}

If $G$ is the Lie group of all maps $f : \mathbb{R}^d\rightarrow
\mathbb{R}^d$ where $f(x) = Ax + b$ for some diagonal matrix $A$
with positive entries and $b\in \mathbb{R}^d$ we obtain a
generalization of Theorem~\ref{thm:prva} to the case of $G$-orbits
of convex bodies, including for example the case of axes-aligned
ellipsoids.

\begin{rem}\label{rem:justify-cubes}{\rm
Perhaps as a justification of treating separately the cases of
cuboids (Theorems~\ref{thm:prva} and \ref{thm:druga}) from the
general case (Theorem~\ref{thm:Lie-group}), here we argue that
after all there is something special about the family
$\mathcal{P}_d$.

Suppose that $f_i : \mathbb{R}\rightarrow \mathbb{R} \,
(i=1,\ldots, d)$ is a family of homeomorphisms. If $f = \prod f_i$
is the associated auto-homeomorphism of $\mathbb{R}^d$ then $f$
obviously sends cuboids to cuboids. Similarly if $f_i :
\mathbb{R}\rightarrow (0,1)$  are homeomorphisms, the associated
product homeomorphism $f : \mathbb{R}^d \rightarrow (0,1)^d$ sends
bijectively the cuboids from $\mathbb{R}^d$ to cuboids from
$(0,1)^d$. By restricting on $(0,1)^d$ the functions described in
the Example~\ref{exam:treci-primer}, normalizing and pulling back
to $\mathbb{R}^d$ by $f$, we can easily construct $2d$ probability
measures on $\mathbb{R}^d$ which distinguish cuboids one from
another. As a consequence we can prove that $\nu(\mathcal{P}_d,
\mathcal{P})= 2d$ (see Remark~\ref{rem:general-nu}) where
$\mathcal{P}$ is the family of probability measures on
$\mathbb{R}^d$. This is a result that does not obviously hold for
other classes $\mathcal{F}$ covered by
Theorem~\ref{thm:Lie-group}. }
\end{rem}

\section{Some exact values for $\nu_{rel}(\mathcal{F})$}

In this section we show that the invariant
$\nu_{rel}(\mathcal{F})$ can be evaluated for many other classes
of convex bodies (measurable sets) in $\mathbb{R}^d$. In
particular we show that the same bounds that we determined in the
case of cuboids (Theorem~\ref{thm:prva}), also hold in the case of
ellipsoids and other centrally symmetric, axes-aligned bodies in
$\mathbb{R}^d$.

\medskip
We begin with some preliminary definitions. Let $\chi :
\mathbb{R}^d\rightarrow \mathbb{R}$ be a non-negative, bounded,
measurable function with bounded support in $\mathbb{R}^d$. Our
main example of such a function is the indicator function $\chi_K$
of a convex body $K\subset \mathbb{R}^d$.

We put some emphasis on the case of centrally symmetric convex
bodies. Motivated by that we say that $\chi :
\mathbb{R}^d\rightarrow \mathbb{R}$ is an even function, relative
$a\in \mathbb{R}^d$, if
\begin{equation}\label{eqn:chi}
\chi(a+x) = \chi(a-x) \quad \mbox{for each }\, x\in \mathbb{R}^d.
\end{equation}
Let $G_{DT}$ be the group of all transformations of $\mathbb{R}^d$
generated by translations and positive axes-aligned dilatations.
More explicitly, $L\in G_{DT}$ if there exists a diagonal matrix
$A = {\rm diag}\{C_1,\ldots, C_d\}$ with positive entries, and a
vector $b\in \mathbb{R}^d$ such that $L(x) = A(x) + b$ for each
$x\in \mathbb{R}^d$.

\medskip
It is not difficult to check that the geometric dimension of the
family $\mathcal{F}_K = \{L(K)\mid L\in G_{DT}\}$ is equal to
$2d$. In light of Theorem~\ref{thm:Lie-group} we know that
$\nu_{rel}(\mathcal{F}_K)\geq 2d$. The following proposition
establishes the opposite inequality in the case of centrally
symmetric convex bodies.

\begin{theo}\label{thm:opposite}
Let $K$ be a centrally symmetric convex body in $\mathbb{R}^d$.
Let $$\mathcal{F}_K = \{L(K)\mid L\in G_{DT}\}$$ be the associated
family of all convex bodies obtained from $K$ by successive
positive, axes-aligned dilatations and translations. Then
$\nu_{rel}(\mathcal{F}_K)\leq 2d$. More explicitly, the required
relaxed $(2d)$-coloring of $\mathbb{R}^d$ is provided by the
measures $d\mu_i = \phi_i dm$ with the following density
functions,
\begin{equation}\label{eqn:density}
\phi_0 = 1, \quad \phi_i(x) = x_i \quad (i=1,\ldots, d),\quad
\phi_{d+i} = x_i^2 \quad (i=1,\ldots, d-1).
\end{equation}
\end{theo}

Before we commence the proof of Theorem~\ref{thm:opposite} we
establish the following lemma.

\begin{lema}\label{lema:moment}
Suppose that $\alpha : \mathbb{R}\rightarrow \mathbb{R}$ is a
non-negative, integrable function such that for some $r\in
\mathbb{R},\, \alpha(r-x) = \alpha(r+x)$ for each $x\in
\mathbb{R}$. Assume that $\int_\mathbb{R}\alpha > 0$. Suppose that
$\mathfrak{u}_1(x) = a_1x+b_1$ and $\mathfrak{u}_2(x) = a_2x+b_2$
are two increasing linear functions $(a_1, a_2 > 0)$ such that,
\begin{equation}\label{eqn:moment-lema}
\int_{\mathbb{R}}\mathfrak{u}_1(x)\alpha(x)\, dx =
\int_{\mathbb{R}}\mathfrak{u}_2(x)\alpha(x)\, dx \quad \mbox{and}
\quad \int_{\mathbb{R}}\mathfrak{u}_1^2(x)\alpha(x)\, dx =
\int_{\mathbb{R}}\mathfrak{u}_2^2(x)\alpha(x)\, dx.
\end{equation}
Then $(a_1,b_1)=(a_2,b_2)$, i.e.\ $\mathfrak{u}_1(x) =
\mathfrak{u}_2(x)$ for each $x\in \mathbb{R}$.
\end{lema}

\medskip\noindent
{\bf Proof:} Without loss of generality we may assume that $r=0$
which means that $\alpha$ is an even function. Since
$\mathfrak{u}_i(x) + \mathfrak{u}_i(-x)= 2b_i$, it easily follows
from the first equality in (\ref{eqn:moment-lema}) that $b_1=b_2$.

For contradiction let us assume that $a_2> a_1 >0$. It follows
that $\mathfrak{u}_2^2(x) + \mathfrak{u}_2^2(-x) >
\mathfrak{u}_1^2(x) + \mathfrak{u}_1^2(-x)$ for each $x\neq 0$,
which is essentially a consequence of the convexity of the
function $x\mapsto  |x|^2$. However this is in contradiction with
the second equality in (\ref{eqn:moment-lema}) since $\alpha$ is
non-negative and by a change of variables,
\begin{equation}\label{eqn:moment-lema-2}
\int_{\mathbb{R}}[\mathfrak{u}_1^2(x) +
\mathfrak{u}_1^2(-x)]\alpha(x)\, dx =
\int_{\mathbb{R}}[\mathfrak{u}_2^2(x) +
\mathfrak{u}_2^2(-x)]\alpha(x)\, dx
\end{equation}

\medskip\noindent
{\bf Proof of Theorem~\ref{thm:opposite}:} For a given $\chi :
\mathbb{R}^d \rightarrow \mathbb{R}$ let $\mathcal{F}_\chi =
\{\chi\circ L \mid L \in G_{DT}\}$. For example if $\chi = \chi_K$
is the indicator function of a convex body $K$ then
$\mathcal{F}_{\chi_K}$ is the set of indicator functions of
elements from $\mathcal{F}_K$.

We prove a slightly more general statement by showing that if
$\chi : \mathbb{R}^d\rightarrow \mathbb{R}$ is an even,
non-negative, integrable function with non-zero integral, then
each two elements $\chi_1, \chi_2 \in \mathcal{F}_\chi$ are
measure discernible in the sense that,
\begin{equation}\label{eqn:chi-chi}
\int_{\mathbb{R}^d}\phi_i\chi_1 = \int_{\mathbb{R}^d}\phi_i\chi_1
\quad \mbox{\rm for each } i\in\{0,1,\dots, 2d-1\}\quad
\Rightarrow \quad \chi_1 = \chi_2.
\end{equation}

Suppose that $\chi_i(x) = \chi(L_i^{-1})(x) \, (i=1,2)$ where
$L_i(x) = A_i(x) + b_i$ and
$$
A_1 = {\rm diag}(c_1',\ldots, c_d'),\, A_2 = {\rm
diag}(c_1'',\ldots, c_d''),\, b_1 = (b_1',\ldots, b_d'),\, b_2 =
(b_1'',\ldots, b_d'').
$$
The assumption from (\ref{eqn:chi-chi}) on the functions $\chi_1$
and $\chi_2$ is by a change of variables equivalent to,
\begin{equation}\label{eqn:products}
\prod_{i=1}^d c_i' \int_{\mathbb{R}^d} \phi_i(L_1(x))\chi(x)\, dx
= \prod_{i=1}^d c_i'' \int_{\mathbb{R}^d} \phi_i(L_2(x))\chi(x)\,
dx.
\end{equation}
Since $\phi_0=1$, and by assumption $\int \chi \neq 0$, we deduce
from (\ref{eqn:products}) that $\prod_{i=1}^d c_i' = \prod_{i=1}^d
c_i''$.

The functions $\phi_j(x) = x_j$ and $\phi_{d+j}(x) = x_j^2$ depend
only on the variable $x_j$. If $i\in\{j, d+j\}$ then the equality
(\ref{eqn:products}) (after cancelling out the products) can be
rewritten as follows,
\begin{equation}\label{eqn:products-2}
\int_{\mathbb{R}} \phi_i(L_1(x))\widehat{\chi}(x)\, dx =
\int_{\mathbb{R}} \phi_i(L_2(x))\widehat{\chi}(x)\, dx
\end{equation}
where $\widehat{\chi}$ is the density of the measure on
$\mathbb{R}$ obtained as the pushforward of the measure $\chi\,
dx$ (defined on $\mathbb{R}^d$) with respect to the projection on
the $x_j$-axis.

Since $\phi_j(L_1(x)) = c_j'x_j + b_j' $ and $\phi_{d+j}(L_1(x)) =
(c_j'x_j + b_j')^2$ the equalities (\ref{eqn:products-2}) provide
exactly the input needed for the application of
Lemma~\ref{lema:moment}. As a consequence we have the equality
$(c_j', b_j') = (c_j'', b_j'')$ for each $j=1,2,\ldots, d-1$.

From here and the equality of products we observe that $c_d' =
c_d''$. By choosing the remaining unused function $\phi_d(x) =
x_d$, and by one more application of (\ref{eqn:products-2}), we
deduce that $b_d'=b_d''$ which completes the proof of the
proposition. \hfill $\square$

\begin{cor}\label{cor:mu-konv-tela}
Suppose that $K$ is a centrally symmetric convex body in
$\mathbb{R}^d$. If $\mathcal{F}_K = \{L(K)\mid L\in G_{DT}\}$ is
the associated family of all convex bodies obtained by successive
dilatations (positive, axis-aligned) and translations, then
\begin{equation}
\nu_{rel}(\mathcal{F}_K) = 2d.
\end{equation}
\end{cor}

\end{document}